\def \latest {2018--08--17}

\documentclass[12pt]{article}
\usepackage {graphics}
\usepackage{diagbox,multirow}
\usepackage {pstricks,pst-node,pst-coil}
\usepackage{epsfig}
\usepackage{amssymb,amsfonts,amsmath}%
\usepackage{latexsym}

\def\zet{\mathbb{Z}}
\def\vi{V_4}
\def \vcor {$\vi$-cordial}

\def\fk2{\lfloor\frac{k}{2}\rfloor}
\def\ck2{\lceil\frac{k}{2}\rceil}
\def\k1{\lfloor\frac{k+1}{2}\rfloor}

\newcommand{\qed}{\hfill \rule{.1in}{.1in}}

\newtheorem{theorem}{Theorem}
\newtheorem{conjecture}[theorem]{Conjecture}
\newtheorem{proposition}[theorem]{Proposition}
\newtheorem{problem}[theorem]{Problem}

\begin{document}
\title{$\zet_2\times\zet_2$-cordial cycle-free hypergraphs}
\author{Sylwia Cichacz$^{1}$, Agnieszka G\"{o}rlich$^1$, Zsolt
Tuza$^{2,3}$\\
$^1$AGH University of Science and Technology,\\ 
al. A. Mickiewicza 30, 30-059 Krakow, Poland\\ 
$^2$Alfr\'{e}d R\'{e}nyi Institute of Mathematics, Hungarian Academy of Sciences,\\
H-1053 Budapest, Re\'{a}ltanoda u.~13--15;\\
$^3$Department of Computer Science and Systems Technology,\\
University of Pannonia, H-8200 Veszpr\'{e}m, Egyetem
u.~10, Hungary}

\date{\small Latest update on \latest }
 \maketitle

\begin{abstract}
Hovey introduced $A$-cordial labelings  as a  generalization of cordial and harmonious labelings  \cite{Hovey}. If $A$ is an Abelian group, then a labeling $f \colon V (G) \rightarrow A$ of the vertices of some graph $G$ induces an edge labeling on $G$; the edge $uv$ receives the label $f (u) + f (v)$. A graph $G$ is $A$-cordial if there is a vertex-labeling such that (1) the vertex label classes
differ in size by at most one and (2) the induced edge label classes differ in size by at most one.

The problem of $A$-cordial labelings of graphs can be naturally extended for hypergraphs. It was shown that not every $2$-uniform hypertree (i.e., tree) admits a $\zet_2\times\zet_2$-cordial labeling \cite{Pechnik}. The situation changes if we consider $p$-uniform hypetrees for a bigger $p$. We prove that a $p$-uniform hypertree is $\zet_2\times\zet_2$-cordial for any $p>2$, and so is every path hypergraph in which all edges have size at least~3. The property is not valid universally in the class of hypergraphs of
 maximum degree~1, for which we provide a necessary and sufficient condition.

\end{abstract}


\section{Introduction}

A {\it hypergraph} $H$ is a pair $H = (V,E)$ where $V$ is a set of
vertices and $E$ is a set of non-empty subsets of $V$ called
hyperedges. The order
 (number of vertices) of a hypergraph $H$ is denoted by $|H|$ and
the size (number of edges)
 is denoted by $\|H\|$. If all edges have the same
cardinality $p$, the hypergraph is said to be {\it $p$-uniform}. Hence a
graph is $2$-uniform hypergraph. The degree of a vertex $v$, denoted
by $d(v)$, is defined as $d(v) = |\{e \in E : v\in e\}|$; i.e., the
degree of $v$ is the number of edges to which it belongs. Two
vertices in a
hypergraph are adjacent if there is an edge containing both of them.

In order to avoid some trivialities, we assume
 in most of this paper that every edge of a hypergraph has at least two vertices.
The only exception will be Section~\ref{s:match}.

A {\it walk} in a hypergraph is a sequence $v_0, e_1, v_1,
\ldots, v_{n-1}, e_n, v_n$, where $v_i \in  V$, $e_i \in E$ and
$v_{i-1},v_i \in e_i$ for all $i$. We define a {\it path} in a
hypergraph to be a walk with all $v_i$ distinct and all $e_i$
distinct. A {\it cycle} is a walk containing at least two edges, all
$e_i$ are distinct and all $v_i$ are distinct except $v_0 = v_n$. A~hypergraph is connected if for every pair of its vertices $v,u$,
there is a path starting at $v$ and ending at $u$.
A {\it hypertree} is a connected hypergraph with no cycles.

A \textit{star} is a hypertree in which one vertex --- called the
 \textit{center} of the star --- is contained in
 all edges (and the edes are mutually disjoint outside this vertex).
Observe that a $p$-uniform hypertree with $\|T\|$ edges
always has exactly $1 + (p - 1) \|T\|$ vertices.
An even simpler structure is a \textit{matching} --- frequently
 called `packing' in the literature --- in which
 any two edges are vertex-disjoint.
(Here we allow that isolated vertices may also occur.)

\bigskip

For a $p$-uniform hypergraph $H=(V,E)$,  an Abelian group $A$ and an $A$-labeling
$c: V \to A$ let $v_c(a)=|c^{-1}(a)|$.
  The labeling $c$
is said to be {\it $A$-friendly} if $|v_c(a)-v_c(b)| \leq 1$ for any
$a,b \in A$. The labeling $c$ induces an edge
labeling $c^*:E \to A$  defined by $c^*(e)=\sum_{v \in
e}c(v)$.
 Let $e_{c^*}(a)=|{c^*}^{-1}(a)|$. A hypergraph
is said to be {\it $A$-cordial} if it admits an $A$-friendly
labeling $c$ such that  $|e_{c^*}(a)-e_{c^*}(b)|\leq 1$ for any
$a,b \in A$. Then we say that the edge labeling $c^*$ is
{\it $A$-cordial}.

Cordial labeling of graphs was introduced by Cahit \cite{Cahit} as a weakened
version of graceful labeling and harmonious labeling. This notion was generalized by Hovey for any Abelian group of order $k$ \cite{Hovey}.
So far research on $A$-cordiality has mostly focused on the case where $A$ is cyclic  and so called $k$-\textit{cordial}.
Hovey \cite{Hovey} showed that all caterpillars are $k$-cordial for all $k$ and all trees are $k$-cordial for $k=3,4,5$. Moreover he showed that cycles are $k$-cordial for any odd $k$. He raised the conjectures
 that if $H$ is a tree graph, it is $k$-cordial for every $k$, and that all connected graphs are $3$-cordial \cite{Hovey}. In the last twenty-five years there was little progress towards a solution to either of these
conjectures. However, Driscoll, Krop and Nguyen proved recently that all trees are 6-cordial \cite{ref_DriKroNgu}.

Note that this result  does not extend even to the smallest noncyclic group,
the Klein four-group
(i.e., $\vi=\zet_2 \times \zet_2$); the paths $P_4$ and $P_5$ are not
$\vi$-cordial what is shown in the following theorem.

\begin{theorem}[\cite{Pechnik}]\label{path}
The path $P_n$ is $\vi$-cordial unless $n \notin \{4, 5\}$.
\end{theorem}

In \cite{CGT} we investigated a problem analogous to Hovey's problem for hypertrees (connected hypergraphs without cycles) and presented various sufficient conditions on $H$ to be $k$-cordial. From our theorems it
follows that  every uniform hyperpath is $k$-cordial
for any $k$, and  every $k$-uniform hypertree is $k$-cordial.  We  conjectured that all  hypertrees  are $k$-cordial for all $k$.   Recently Tuczy\'{n}ski, Wenus and W\c{e}sek  proved this conjecture for $k=2,3$ \cite{Tuczu}.

However, a $2$-uniform hypertree is not $\vi$-cordial in general by Theorem~\ref{path}.

In this paper we show that such counterexamples no longer exist in case of  $p$-uniform hypertrees for $p\geq3$.
Namely, we prove that any $p$-uniform hypertree is $\vi$-cordial for all $p\geq3$. Beyond that, for stars we can even drop the condition of uniformity.
We also characterize $\vi$-cordial hypergraphs whose edges are mutually disjoint
 (i.e., matchings).


\section{ Extension Lemma and uniform hypertrees
}

We begin this section with some sufficient conditions under which a
 $\vi$-cordial labeling can be derived from that of a subhypergraph.
This result will be applied later in several situations, leading to substantial
 shortening of various arguments.
We use it first for uniform hypertrees, proving that all of them are $\vi$-cordial.

Before we present the results, we introduce a notation
for convenience. Let the edge set of the hypergraph under consideration be
 $E = \{e_1, e_2, \ldots,e_m\}$.
  For all $1 \le i \leq m$, let us denote $X_{i} = \bigcup_{1\leq j \leq i} e_j$.
We will assume without loss of generality that the edges are indexed in such a way
 that $e_i$ meets at most one connected component of the subhypergraph with
 vertex set $X_{i-1}$ and edge set $\{e_1, \ldots,e_{i-1}\}$.
In particular, for hypertrees it means that each $e_i$ has exactly one vertex
 in common with the set $X_{i-1}$; hence every $\{e_1, e_2,\ldots , e_i \}$
 forms a hypertree in which $e_i$ is a pendant edge.
For hypertrees it can also be assumed that $e_m$ is the last edge
 in a longest path in $T$.

%

\begin{theorem}     \label{l:ext}
{\bf (Extension Lemma)} \
Let $H=(V,E)$ be a hypergraph with edge set $E=\{e_1,\dots,e_m\}$,
 and let $e_m^- := e_m \setminus  (e_1\cup\cdots\cup e_{m-1})$.
  Assume that $|e_m^-|\ge 2$, and that the following conditions hold:
 \begin{enumerate}
   \item If $|V|\equiv0 \pmod 4$, then $m\equiv1 \pmod 4$.
   \item If $|V|\equiv2 \pmod 4$, then $m\not\equiv0 \pmod 4$.
   \item If $|V|\equiv3 \pmod 4$ and $|e_m^-|=2$, then $m\not\equiv0 \pmod 4$.
 \end{enumerate}
 If the hypergraph $H^-$ obtained from $H$ by omitting $e_m$ from $E$
  and deleting the vertices of $e_m^-$ from $V$ is $\vi$-cordial,
   then $H$ is $\vi$-cordial.
\end{theorem}
\textit{Proof.}\\
Assume that $c'$ is a $\vi$-labeling of $X_{m-1}$ that induces a $\vi$-cordial
labeling $c'^*$ of $H^-$. If $m-1\equiv 0 \pmod 4$,
then every $\vi$-friendly extension of $c'$ to $X_m=V$ verifies that
$H$ is $\vi$-cordial. Otherwise, if $m\not \equiv 1 \pmod 4$,
assumption \textit{1.}\ of the theorem implies $|X_m|\not \equiv 0 \pmod 4$. If $|X_{m-1}| \not \equiv 0 \pmod 4$, we first assign $a:=4-(|X_{m-1}|
\pmod 4)$ vertices of $e_m^-=e_m\setminus X_{m-1}$ to those
elements of $\vi$ which occur on one fewer vertices of
$X_{m-1}$ than the other $4-a$ elements. Here 
 $1\le a\le 3$,
  and the step is feasible unless $|e_m^-|=2$ and $a=3$,
 because apart from this exception
$|e_m^-|\ge a$ holds and there is enough room to have the current
partial labeling completely balanced for the elements of
$\vi$.

Suppose first that either $|e_m^-|\ge 3$ or $a\le 2$.
Let $b=|e_m^-|-a$ denote the number of vertices
unlabeled so far. We next distribute equally the elements of
$\vi$ on $b - (b\pmod4)$ vertices of $e_m$. There
still remain some $r$ unlabeled vertices in $e_m$, where $1\le
r\le 3$ since $|X_m| \not \equiv 0 \pmod 4$. We choose
$q\in\vi$ such that the current partial sum on $e_m$ plus
$q$  occurs fewer times than some other label(s) in
$c'^*$ on the edge set $e_1,\dots,e_{m-1}$.
  By assumption \textit{2.}\ that $|X_m|\equiv2 \pmod 4$ implies
   $m\not\equiv0 \pmod 4$, we can take $q\neq(0,0)$ if $r=2$.
Therefore we can easily select $r$ distinct elements $l_1,\dots,l_r\in \vi$
such that $l_1+\dots+l_r = q$. Assigning them to
the remaining vertices, a $\vi$-cordial labeling of the entire $T$
is obtained.

Consider now the case $|e_m^-|=2$ with $a=3$.
Here $a=3$ means that $|X_{m-1}| \equiv 1 \pmod 4$,
 and then $|e_m^-|=2$ yields $|X_m| \equiv 3 \pmod 4$. Hence so far
 three elements of $\vi$ are used one fewer than the fourth element,
 and we have to use two of them on the unlabeled vertices of $e_m$.
Now $m\not\equiv0 \pmod 4$ by assumption \textit{3.}, thus at least two
 sums are feasible on $e_m$.
Consequently, by the pigeonhole principle, one of two feasible sums
 coincides with one of three sums which can be generated by the
 sum of labels on the vertices in $X_{m-1}\cap e_m$ together with the
  pairs of the three usable elements of $\vi$.
~\qed

\begin{theorem}\label{uniform}

Let $p\ge 3$. Then every $p$-uniform hypertree is $\vi$-cordial.
\end{theorem}
\textit{Proof.}\\
 The theorem obviously holds for any hypertree with size one,
this case is the anchor of induction.
  Let $T$ be a $p$-uniform hypertree with size $m=\|T\|\ge 2$ and assume that the theorem holds for every $p$-uniform hypertree with size less than $m$.  Let $T'=T-\{e_m\}$ be the $p$-uniform hypertree with vertex set $V'=X_{m-1}$. By induction there exists a $\vi$-friendly labeling $c'$ for $T'$ which induces a $\vi$-cordial labeling ${c'}^*$. Below we show that $c'$ can be extended to a $\vi$-friendly labeling $c$ of $T$ in such a way that $c$ induces $\vi$-cordial labeling for $T$.

Recall that we have $|T|=(p-1)\|T\|+1$, therefore the residue of $|X_m|$ modulo~4
 is obtained according to Table \ref{tab:1}.
Column $m\equiv 0$ shows that the second and third conditions in
 Theorem \ref{l:ext} automatically hold, moreover only one of the two
 occurrences of 0 violates the first condition.
Hence, to complete the proof, we may restrict our attention to
 $p\equiv2 \pmod 4$ and $m\equiv 3\pmod 4$,
 in which case we have $|X_m|\equiv 0\pmod 4$.
We will consider three subcases.

\renewcommand{\arraystretch}{1.5}
\begin{table}[htp]
	\begin{center}%
\begin{tabular}
[c]{|r||rrrr|}\hline
(mod 4)~ & ~$m\equiv ~~~ 0$ & 1 & 2 & 3~ \\ \hline\hline
~$p\equiv ~~~ 0$~ & 1 & 0 & 3 & 2~ \\ 
1~ & 1 & 1 & 1 & 1~ \\
2~ & 1 & ~2 & ~3 & ~0~ \\ 
3~ & 1 & 3 & 1 & 3~ \\ 
 \hline
\end{tabular}
\end{center}
\caption{The value of $|X_m| \pmod 4$.}\label{tab:1}
\end{table}

\medskip


\noindent
{\it Case 1. } $e_{m-1}\cap e_m \neq \emptyset$

\noindent
Note that in this situation $T''=T-\{e_{m-1},e_{m}\}$ is a $p$-uniform hypertree with the vertex set $V''=X_{m-2}$.
By the induction hypothesis there exists a $\vi$-friendly labeling $c''$ for $T''$ which induces a $\vi$-cordial labeling ${c''}^*$.
We show that $c''$ can be extended to a $\vi$-friendly labeling $c$ of $T$ in such a way that $c$ induces $\vi$-cordial labeling for $T$.
Note that in this case there are exactly two elements $x,y \in \vi$ that occur  one time fewer in the labeling $c''$ of the vertices of $T''$  than the other two elements of $\vi$; and there is exactly one element $z\in \vi$ that occurs one time more in the labeling of the edges of $T''$ induced by $c''$ than the other three elements of $\vi$.
 Let $e_{m-1}=\{v,v_1^{m-1},v_2^{m-1},\ldots,v_{p-1}^{m-1}\}$ and  $e_m=\{v,v_1^m,v_2^m,\ldots,v_{p-1}^m\}$.

Suppose first that $X_{m-2}\cap e_{m-1}=\{v\}$.
If now  $z\not \in \{x+c''(v), y +c''(v)\}$ then we put label $x$ on $v_{1}^{m-1}$ and $y$ on $v_1^m$, and on the remaining vertices of the edges $e_{m-1}$ and $e_m$  each element of $\vi$ exactly $(p-2)/4$ times.
Obviously we obtain a $\vi$-cordial labeling of $T$.
If $z \in \{x+c''(v), y +c''(v)\}$, then there exists $\alpha \in \vi$ such that  $z \not\in \{x+c''(v)+\alpha, y +c''(v)+\alpha\}$.
Label vertices as follows: $v_{1}^{m-1}$ by $x$,  $v_{2}^{m-1}$ by $\alpha$, and $v_{3}^{m-1}$, $v_{4}^{m-1}$, $v_{5}^{m-1}$ by the elements $(0,1),(1,0),(1,1)$, whereas $v_{1}^m$ by $y$,  $v_{2}^m$ by $(0,0)$, and $v_{3}^m$, $v_{4}^m$, $v_{5}^m$ by the elements of $\vi-\{\alpha\}$; and on the remaining vertices put each element of $\vi$ exactly $(p-6)/4$ times in each of $e_{m-1}$ and $e_m$.

Suppose now that $X_{m-2}\cap e_{m-1}\neq \{v\}$,
say  $X_{m-2}\cap e_{m-1}= \{v_1^{m-1}\}$.
We can assume that $x+c''(v_1^{m-1})\neq z$ because $y\neq x$.
Label 
  $v_{2}^{m-1}$ by $x$ and put on the remaining vertices of the edge $e_{m-1}$   each element of $\vi$ exactly $(p-2)/4$ times in such a way that $y+c(v)\not \in \{z,x+c''(v_1^{m-1})\}$.
  Label now $v_{1}^{m}$ by $y$ and put on the remaining vertices of the edge $e_{m}$   each element of $\vi$ exactly $(p-2)/4$ times.

\medskip


\noindent
{\it Case 2. } $e_{m-1}\cap e_m = \emptyset$

\noindent
One can easily see (and it also follows from the inductive step described below)
 that if $m=3$, then the hypertree (path) $T$ is $\vi$-cordial.
 Therefore we can assume that $m\geq7$.
Observe that this time $T''=T-\{e_{m-2},e_{m-1},e_{m}\}$ is a $p$-uniform hypertree with the vertex set $V''=X_{m-3}$.
By induction there exists a $\vi$-friendly labeling $c''$ for $T''$ which induces a $\vi$-cordial labeling ${c''}^*$.
Note that in this case there are exactly three elements $x,y,z \in \vi$ that occur  one time fewer in the labeling $c''$ of vertices $T''$  than the other element of $\vi$, and  all the elements of $\vi$ occur the same times  in the labeling of edges of $T''$ induced by $c''$.
We show that the labeling $c''$ can be extended to a $\vi$-friendly labeling $c$ of $T$ in such a way that $c$ induces a $\vi$-cordial labeling for $T$.

Assume first that
 $e_{m-2}\cap e_m \neq \emptyset$ and $e_{m-2}\cap e_{m-1} \neq \emptyset$.
Let us denote $e_{m-2}=\{v_1^{m-2},v_2^{m-2},\ldots,v_{p}^{m-2}\}$, $e_{m-1}=\{v_2^{m-2},v_1^{m-1},v_2^{m-1},\ldots,v_{p-1}^{m-1}\}$ and  $e_m=\{v_3^{m-2}, v_1^m,v_2^m,\ldots,v_{p-1}^m\}$ such that $X_{m-3}\cap e_{m-2}=\{v_1^{m-2}\}$;
for the moment we assume that $v_1^{m-2}\not\in\{v_2^{m-2},v_3^{m-2}\}$.
Put label $x$ on the vertex $v_2^{m-2}$,
 and on the remaining vertices of the edge $e_{m-2}$  each element of $\vi$ exactly $(p-2)/4$ times in such a way that $c(v_2^{m-2})=c(v_3^{m-2})$. For the edges $e_{m-1}$ and $e_m$ proceed the same way now as in Case 1.

 In the other situation, if $X_{m-3}\cap e_{m-2}$ coincides with
 $e_{m-2}\cap e_{m-1}$, we apply essentially the same strategy,
 imposing the condition that the vertex $e_{m-2}\cap e_m$
 gets the label $c''(v_1^{m-2})$.

Next, let $e_{m-2}\cap e_m =\emptyset$ and $e_{m-2}\cap e_{m-1} \ne \emptyset$.
This situation can be reduced to Case 1 by a modification of the indexing of the edges, viewing $e_{m-1}$ as the new $e_m$, also $e_{m-2}$ as the new $e_{m-1}$, and the old $e_m$ (which is disjoint from both other edges) as the new $e_{m-2}$.
Using the new indices we have $e_{m-1}\cap e_m \ne \emptyset$, which has already been settled.
A similar re-indexing works if
  $e_{m-2}\cap e_{m-1} =\emptyset$ and $e_{m-2}\cap e_m \ne \emptyset$.

Finally, assume that $e_{m-2}\cap e_m =\emptyset$ and $e_{m-2}\cap e_{m-1} = \emptyset$. Then let   $e_{m-2}=\{v_1^{m-2},v_2^{m-2},\ldots,v_{p}^{m-2}\}$, $e_{m-1}=\{v_1^{m-1},v_2^{m-1},\ldots,v_{p}^{m-1}\}$ and  $e_m=\{v_1^m,v_2^m,\ldots,v_{p}^m\}$ such that $X_{m-3}\cap e_{m-2}=\{v_1^{m-2}\}$, $X_{m-3}\cap e_{m-1}=\{v_1^{m-1}\}$ and $X_{m-3}\cap e_{m}=\{v_1^{m}\}$. Suppose first that $|\{c''(v_1^{m-2}), c''(v_1^{m-1}),c''(v_1^{m})\}|<3$, then without loss of generality we can assume that $c''(v_1^{m-1})=c''(v_1^{m})$. Put label $x$ on the vertex $v_{p}^{m-2}$  and on the remaining vertices of the edge $e_{m-2}$  each element of $\vi$ exactly $(p-2)/4$ times. For the edges $e_{m-1}$ and $e_m$ proceed the same way now as in Case 1.

Otherwise, if $\{c''(v_1^{m-2}), c''(v_1^{m-1}),c''(v_1^{m})\}=\{a,b,c\}$
 is a set of three distinct labels, we let $\beta=\vi-\{a,b,c\}$.
On $p-2$ vertices in each of $e_{m-2},e_{m-1},e_m$ we distribute the elements
 of $\vi$ equally, using $(p-2)/4$ times each.
The current partial sums on these edges are $a,b,c$, and we need to assign
 $x,y,z$ (one of them in each edge) in a way that the sums remain mutually distinct.
If $\beta\notin\{x,y,z\}$, then in fact $\{a,b,c\}=\{x,y,z\}$, and we can
 obviously create the sums $x+y$, $y+z$, and $z+x$, which satisfy the conditions.
Else, if say $\beta=x$, we have $\{a,b,c\}=\{a,y,z\}$ where $a\ne x$.
We then create two nonzero sums $a+y$ and $y+x$, and the zero sum $z+z$.
The corresponding labeling satisfies the conditions and completes the proof
 of the theorem.~\qed

\section{Stars, matchings, paths}

In this section we consider hypergraphs also with smaller edges than
 in the previous sections, because even such extensions allow
 characterizations for the existence of $\vi$-cordial labelings
 in some subclasses.
In particular, stars need no restriction, whereas $\vi$-cordial hypergraphs of
 maximum degree 1 admit a simple characterization.
The case of paths seems to be more complicated to handle, here we only exhibit an
 infinite family which is not $\vi$-cordial.


\subsection{Stars}

Recall that the edge set of a star is a collection of sets of size at least 2
 each, which are mutually disjoint apart from a single vertex which is
 contained in all of them.
Hence each edge $e_i$ contains precisely $|e_i|-1$ private vertices,
 and with the notation of the Extension Lemma (Theorem \ref{l:ext})
 we have $|e_m^-|=|e_m|-1$, no matter which indexing order $e_1,\dots,e_m$
  of the edges we take.

\begin{theorem}   \label{t:star}
 Every star is $\vi$-cordial.
\end{theorem}
\textit{Proof.}\\
Let $H$ be a star with $m$ edges $e_1,\dots,e_m$.
We can associate the $m$-tuple $(f_1,\dots,f_m)$ of integers with $H$,
 where $f_i=|e_i|-1$ for all $1\le i\le m$.
It is clear that every $m$-tuple of positive integers uniquely determines
 the corresponding star up to isomorphism, moreover $|H|=1+\sum_{i=1}^m f_i$.
This representation can further be simplified to one which still determines
 $H$, namely we can denote by $m_k$ the number of indices $i$ such that
 $f_i=k$.

The proof will be an induction on $|H|$, anchored by approximately 30
 small cases.
We are going to introduce several reductions, along which it will turn out
 which of the small cases are relevant to be checked separately.
Below we describe situations and explain why they are reducible.
  \begin{itemize}
    \item[(1)] If there is a $k\ge 5$ with $m_k>0$, then it reduces to
      $m_k := m_k -1$ and $m_{k-4} := m_{k-4} +1$.
  \end{itemize}
The reason is that inside an edge with 5 or more non-center vertices
 we can assign four to the elements of $\vi$, hence creating a partial
 sum equal to zero and decreasing $|H|$ by four, still having a star
 with $m$ edges.
Hence it suffices to consider stars represented by 4-tuples
 $(m_1,m_2,m_3,m_4)$.
  \begin{itemize}
    \item[(2)] If there is a $k\le 4$ with $m_k\ge 4$, then it reduces to
      $m_k := m_k -4$.
  \end{itemize}
Assume that $|e_1|=|e_2|=|e_3|=|e_4|=k+1$. Table \ref{tab:2} shows how the
 non-center vertices of $e_1,e_2,e_3,e_4$ can be labeled to induce four
 distinct edge labels, and hence eliminate those four edges.
In this way all remaining stars to be considered are represented by 4-tuples
 $(m_1,m_2,m_3,m_4)\in\{0,1,2,3\}^4$, that is already a finite collection of
 basic configurations.

\renewcommand{\arraystretch}{1.5}
\begin{table}[htp]
{\small
	\begin{center}%
\begin{tabular}
[c]{|c||c|c|c|c|}\hline
 & $k=1$ & $k=2$ & $k=3$ & $k=4$ \\ \hline\hline
$e_1 = (0,0)$ & $(0,0)$ & $(0,0)$, $(0,0)$ & $(0,0)$, $(0,0)$, $(0,0)$ & $(0,0)$, $(0,1)$, $(1,0)$, $(1,1)$ \\ 
 $e_2 = (0,1)$ & $(0,1)$ & $(1,0)$, $(1,1)$ & $(0,1)$, $(0,1)$, $(0,1)$ & $(0,0)$, $(0,1)$, $(0,1)$, $(0,1)$ \\
 $e_3 = (1,0)$ & $(1,0)$ & $(0,1)$, $(1,1)$ & $(1,0)$, $(1,0)$, $(1,0)$ & $(0,0)$, $(1,0)$, $(1,0)$, $(1,0)$ \\ 
 $e_4 = (1,1)$ & $(1,1)$ & $(0,1)$, $(1,0)$ & $(1,1)$, $(1,1)$, $(1,1)$ & $(0,0)$, $(1,1)$, $(1,1)$, $(1,1)$ \\ 
 \hline
\end{tabular}
\end{center}
\caption{Eliminating four edges of equal size. The label of center vertex, when different from $(0,0)$, permutes the edge sums indicated in the first column.}\label{tab:2}
}
\end{table}

  \begin{itemize}
    \item[(3)] If $f_1+f_2+f_3+f_4\equiv 0\pmod 4$, then
     $e_1,e_2,e_3,e_4$ can be eliminated.
    More explicitly, if in each position the 4-tuple $(m_1,m_2,m_3,m_4)$ is
     at least as large as one or more of
     $$
       (0,1,2,1), ~~ (0,2,0,2), ~~ (1,0,1,2), ~~ (1,2,1,0), ~~ (2,0,2,0), ~~ (2,1,0,1)
     $$
     then the configuration is reducible.
  \end{itemize}
Indeed, the condition $f_1+f_2+f_3+f_4\equiv 0\pmod 4$
 actually means that $f_1+f_2+f_3+f_4$ equals 8 or 12,
 because 4 and 16 would only occur as $4\times 1$ and $4\times 4$,
 respectively, and these cases have just been settled by (2).
Simple enumeration yields that there are six possible 4-tuples
 $(f_1,f_2,f_3,f_4)$ apart from permutations.
Table \ref{tab:3} exhibits an ad hoc labeling from the many possibilities
 for each of them, showing that all these subconfigurations can be eliminated.
There is a direct one-to-one correspondence between the 4-tuples
 $(m_1,m_2,m_3,m_4)$ and $(f_1,f_2,f_3,f_4)$, for example
 $(m_1,m_2,m_3,m_4)=(1,0,2,1)$ --- the third case listed above ---
 means $f_1=1$, $f_2=3$, $f_3=3$, $f_4=4$.

\renewcommand{\arraystretch}{1.5}
\begin{table}[htp]
{\footnotesize
	\begin{center}%
\begin{tabular}
[c]{|c||c|c|c|c|}\hline
$(f_1,f_2,f_3,f_4)$ & $e_1 = (0,0)$ & $e_2 = (0,1)$ & $e_3 = (1,0)$ & $e_4 = (1,1)$ \\ \hline\hline
 $(1,1,2,4)$ & $(0,0)$ & $(0,1)$ & $(0,0)$, $(1,0)$ & $(0,1)$, $(1,0)$, $(1,1)$, $(1,1)$ \\
 $(1,1,3,3)$ & $(0,0)$ & $(0,1)$ & $(0,0)$, $(0,1)$, $(1,1)$ & $(1,0)$, $(1,0)$, $(1,1)$ \\
 $(1,2,2,3)$ & $(0,0)$ & $(0,0)$, $(0,1)$ & $(0,1)$, $(1,1)$ & $(1,0)$, $(1,0)$, $(1,1)$ \\ \hline
 $(1,3,4,4)$ & $(0,0)$ & $(0,1)$, $(0,1)$, $(0,1)$ & $(0,0)$, $(1,0)$, $(1,0)$, $(1,0)$ & $(0,0)$, $(1,1)$, $(1,1)$, $(1,1)$ \\
 $(2,2,4,4)$ & $(0,0)$, $(0,0)$ & $(0,0)$, $(0,1)$ & $(0,1)$, $(1,0)$, $(1,0)$, $(1,1)$ & $(0,1)$, $(1,0)$, $(1,1)$, $(1,1)$ \\
 $(2,3,3,4)$ & $(0,0)$, $(0,0)$ & $(0,1)$, $(0,1)$, $(0,1)$ & $(1,0)$, $(1,0)$, $(1,0)$ & $(0,0)$, $(1,1)$, $(1,1)$, $(1,1)$ \\
 \hline
\end{tabular}
\end{center}
\caption{Eliminating four edges whose total number of non-center vertices is 8 or 12.}\label{tab:3}
}
\end{table}

Since the theorem claims $\vi$-cordiality of stars without any exceptions,
 all the situations described above provide an inductive step
  when they occur as subconfigurations.
It follows that, for an anchor of the induction, a $\vi$-cordial labeling
 has to be presented for only those stars which are not reducible by any
 of (1)--(3).
There are 79 such cases, as listed in Table \ref{tab:4}.
Below we show how they can be handled.

\renewcommand{\arraystretch}{1.5}
\begin{table}[htp]
{
\small
	\begin{center}%
\begin{tabular}
[c]{|c|c|c||c|c|c||c|c|c||c|c|c|}\hline
 $(0,0,0,1)$ &  $1,5$ & {\bf O} &  $(0,1,1,1)$ &  $3,10$ & {\bf F} &  $(1,0,0,3)$ &  $4,14$ & {\bf F} &  $(2,0,0,1)$ &  $3,7$ & {\bf R}   \\
 $(0,0,0,2)$ &  $2,9$ & {\bf T} &  $(0,1,1,2)$ &  $4,14$ & {\bf F} &  $(1,0,1,0)$ &  $2,5$ & {\bf T} &  $(2,0,0,2)$ &  $4,11$ & {\bf F}   \\
 $(0,0,0,3)$ &  $3,13$ & {\bf T} &  $(0,1,1,3)$ &  $5,18$ & {\bf T} &  $(1,0,1,1)$ &  $3,9$ & {\bf T} &  $(2,0,0,3)$ &  $5,15$ & {\bf T}   \\
 $(0,0,1,0)$ &  $1,4$ & {\bf O} &  $(0,1,2,0)$ &  $3,9$ & {\bf T} &  $(1,0,2,0)$ &  $3,8$ & {\bf *} &  $(2,0,1,0)$ &  $3,6$ & {\bf R}   \\
 $(0,0,1,1)$ &  $2,8$ & {\bf F} &  $(0,1,3,0)$ &  $4,12$ & {\bf *} &  $(1,0,2,1)$ &  $4,12$ & {\bf F} &  $(2,0,1,1)$ &  $4,10$ & {\bf F}  \\
 $(0,0,1,2)$ &  $3,12$ & {\bf F} &  $(0,2,0,0)$ &  $2,5$ & {\bf T} &  $(1,0,3,0)$ &  $4,11$ & {\bf R} &  $(2,1,0,0)$ &  $3,5$ & {\bf T}  \\
 $(0,0,1,3)$ &  $4,16$ & {\bf F} &  $(0,2,0,1)$ &  $3,9$ & {\bf T} &  $(1,0,3,1)$ &  $5,15$ & {\bf T} &  $(2,1,1,0)$ &  $4,8$ & {\bf *}  \\
 $(0,0,2,0)$ &  $2,7$ & {\bf R} &  $(0,2,1,0)$ &  $3,8$ & {\bf *} &  $(1,1,0,0)$ &  $2,4$ & {\bf *} &  $(2,2,0,0)$ &  $4,7$ & {\bf *}  \\
 $(0,0,2,1)$ &  $3,11$ & {\bf F} &  $(0,2,1,1)$ &  $4,12$ & {\bf F} &  $(1,1,0,1)$ &  $3,8$ & {\bf F} &  $(2,3,0,0)$ &  $5,9$ & {\bf T}  \\
 $(0,0,2,2)$ &  $4,15$ & {\bf F} &  $(0,2,2,0)$ &  $4,11$ & {\bf R}
                                                                    &  $(1,1,0,2)$ &  $4,12$ & {\bf F} &  $(3,0,0,0)$ &  $3,4$ & {\bf O}  \\
 $(0,0,2,3)$ &  $5,19$ & {\bf T} &  $(0,2,3,0)$ &  $5,14$ & {\bf T} &  $(1,1,0,3)$ &  $5,16$ & {\bf T} &  $(3,0,0,1)$ &  $4,8$ & {\bf *}  \\
 $(0,0,3,0)$ &  $3,10$ & {\bf R} &  $(0,3,0,0)$ &  $3,7$ & {\bf R} &  $(1,1,1,0)$ &  $3,7$ & {\bf R} &  $(3,0,0,2)$ &  $5,12$ & {\bf T}  \\
 $(0,0,3,1)$ &  $4,14$ & {\bf F} &  $(0,3,0,1)$ &  $4,11$ & {\bf F} &  $(1,1,1,1)$ &  $4,11$ & {\bf F} &  $(3,0,0,3)$ &  $6,16$ & {\bf F}  \\
 $(0,0,3,2)$ &  $5,18$ & {\bf T} &  $(0,3,1,0)$ &  $4,10$ & {\bf *} &  $(1,1,2,0)$ &  $4,10$ & {\bf *} &  $(3,0,1,0)$ &  $4,7$ & {\bf R}  \\
 $(0,0,3,3)$ &  $6,22$ & {\bf F} &  $(0,3,1,1)$ &  $5,14$ & {\bf T} &  $(1,1,3,0)$ &  $5,13$ & {\bf T} &  $(3,0,1,1)$ &  $5,11$ & {\bf T}  \\
 $(0,1,0,0)$ &  $1,3$ & {\bf O} &  $(0,3,2,0)$ &  $5,13$ & {\bf T} &  $(1,2,0,0)$ &  $3,6$ & {\bf R} &  $(3,1,0,0)$ &  $4,6$ & {\bf *}  \\
 $(0,1,0,1)$ &  $2,7$ & {\bf F} &  $(0,3,3,0)$ &  $6,16$ & {\bf *} &  $(1,2,0,1)$ &  $4,10$ & {\bf F} &  $(3,1,1,0)$ &  $5,9$ & {\bf T}  \\
 $(0,1,0,2)$ &  $3,11$ & {\bf F} &  $(1,0,0,0)$ &  $1,2$ & {\bf O} &  $(1,3,0,0)$ &  $4,8$ & {\bf *} &  $(3,2,0,0)$ &  $5,8$ & {\bf T}  \\
 $(0,1,0,3)$ &  $4,15$ & {\bf F} &  $(1,0,0,1)$ &  $2,6$ & {\bf R} &  $(1,3,0,1)$ &  $5,12$ & {\bf T} &  $(3,3,0,0)$ &  $6,10$ & {\bf R}
  \\
 $(0,1,1,0)$ &  $2,6$ & {\bf R} &  $(1,0,0,2)$ &  $3,10$ & {\bf F} &  $(2,0,0,0)$ &  $2,3$ & {\bf O} &  &  &  \\
 \hline
\end{tabular}
\end{center}
\caption{The 79 cases of $(m_1,m_2,m_3,m_4)$ which are not excluded by
 (1)--(3), the corresponding pairs $m,n$ (number of edges $m=m_1+m_2+m_3+m_4$,
  number of vertices $n=f_1+f_2+f_3+f_4+1$), and a way how they can be settled.
 The 12 cases marked with {\bf *} need labelings to be constructed
  separately.}\label{tab:4}
}
\end{table}

  \begin{itemize}
    \item[{\bf O}] --- Obvious cases are the stars with just one edge
     ($m_1+m_2+m_3+m_4=1$, the $\vi$-cordial labelings are precisely
      the $\vi$-friendly ones) and the star graphs
     ($m_2=m_3=m_4=0$, a labeling is $\vi$-cordial if and only if
      it is $\vi$-friendly on the set of leaves and also on the
      entire vertex set). There are 6 such cases.
  \end{itemize}

  \begin{itemize}
    \item[{\bf T}] --- Trivial reduction applies for stars with 5 edges
     ($m\equiv 1 \pmod 4$, hence the last edge admits any $\vi$-friendly
      extension from a $\vi$-cordial labeling for the first $m-1$ edges);
     and also for stars of order 5 or 9 or 13 ($n\equiv 1 \pmod 4$, hence
      the last vertex can get an arbitrary label needed for a $\vi$-cordial
      extension from $m-1$ edges to $m$ edges).
     This reduction settles 24 cases.
  \end{itemize}

  \begin{itemize}
    \item[{\bf F}] --- Four vertices can be eliminated if $m_4\ge 1$
     and $m_2+m_3+m_4\ge 2$ (here extension goes from $n-4$ to $n$,
      while $m$ remains unchanged). Indeed,
      inside a 5-element edge we can label three non-center vertices with
      $(0,1),(1,0),(1,1)$ while assigning the label $(0,0)$ to a vertex
      in another edge of size at least 3.
     This reduction settles further 24 cases.
  \end{itemize}

  \begin{itemize}
    \item[{\bf R}] --- Reduction applies by Theorem \ref{l:ext} for stars
     with $n\equiv 2 \pmod 4$ unless $m\equiv 0 \pmod 4$; and also with
     $n\equiv 3 \pmod 4$ except when $m\equiv 0 \pmod 4$ and the star contains
     no edges of 4 or 5 vertices (i.e., $m_3=m_4=0$).
     This reduction settles further 13 cases.
  \end{itemize}

  \begin{itemize}
    \item[{\bf *}] --- There are 12 cases not covered by the previous
     considerations; Table \ref{tab:5} exhibits a $\vi$-cordial labeling
     for each of them.
    Although there are several cases, all are very easy to construct.
  \end{itemize}

Together with this last set of labelings {\bf *}, all cases are exhausted
 and the theorem is proved.
~\qed

\renewcommand{\arraystretch}{1.2}
\begin{table}[htp]
{
  \footnotesize
	\begin{center}%
\begin{tabular}
[c]{|c|c||c|c|c|c|}\hline
$(m_1,m_2,m_3,m_4)$ & $m,n$ & $f_i = 1$ & $f_i = 2$ & $f_i = 3$ & $f_i = 4$ \\ \hline\hline
 $(0,1,3,0)$ & $4,12$ &  & $(0,0)$, $(0,0)$ & $(0,1)$, $(0,1)$, $(0,1)$ &  \\
  &  &  &  & $(1,0)$, $(1,0)$, $(1,0)$ &  \\
  &  &  &  & $(1,1)$, $(1,1)$, $(1,1)$ &  \\ \hline
 $(0,2,1,0)$ & $3,8$ &  & $(0,0)$, $(0,1)$ & $(0,1)$, $(1,0)$, $(1,1)$ &  \\
  &  &  & $(0,0)$, $(1,0)$ &  &  \\ \hline
 $(0,3,1,0)$ & $4,10$ &  & $(0,1)$, $(1,0)$ & $(0,0)$, $(0,0)$, $(0,0)$ &  \\
  &  &  & $(0,1)$, $(1,1)$ &  &  \\
  &  &  & $(1,0)$, $(1,1)$ &  &  \\ \hline
 $(0,3,3,0)$ & $6,16$ &  & $(0,0)$, $(0,0)$ & $(0,1)$, $(0,1)$, $(0,1)$ &  \\
  &  &  & $(0,0)$, $(0,0)$ & $(1,0)$, $(1,0)$, $(1,0)$ &  \\
  &  &  & $(0,1)$, $(1,0)$ & $(1,1)$, $(1,1)$, $(1,1)$ &  \\ \hline
 $(1,0,2,0)$ & $3,8$ & $(0,1)$ &  & $(0,0)$, $(0,1)$, $(1,0)$ &  \\
  &  &  &  & $(1,0)$, $(1,1)$, $(1,1)$ &  \\ \hline
 $(1,1,0,0)$ & $2,4$ & $(0,0)$ & $(0,1)$, $(1,0)$ &  &  \\ \hline
 $(1,1,2,0)$ & $4,10$ & $(0,0)$ & $(0,1)$, $(1,0)$ & $(0,0)$, $(0,0)$, $(0,1)$ &  \\
  &  &  &  & $(1,0)$, $(1,1)$, $(1,1)$ &  \\ \hline
 $(1,3,0,0)$ & $4,8$ & $(0,0)$ & $(0,1)$, $(1,0)$ &  &  \\
  &  &  & $(0,1)$, $(1,1)$ &  &  \\
  &  &  & $(1,0)$, $(1,1)$ &  &  \\ \hline
 $(2,1,1,0)$ & $4,8$ & $(0,1)$ & $(0,0)$, $(1,1)$ & $(0,1)$, $(1,0)$, $(1,1)$ &  \\
  &  & $(1,0)$ &  &  &  \\ \hline
 $(2,2,0,0)$ & $4,7$ & $(0,1)$ & $(0,1)$, $(1,0)$ &  &  \\
  &  & $(1,0)$ & $(1,1)$, $(1,1)$ &  &  \\ \hline
 $(3,0,0,1)$ & $4,8$ & $(0,1)$ &  &  & $(0,0)$, $(0,1)$, $(1,0)$, $(1,1)$ \\
  &  & $(1,0)$ &  &  &  \\
  &  & $(1,1)$ &  &  &  \\ \hline
 $(3,1,0,0)$ & $4,6$ & $(0,0)$ & $(0,0)$, $(1,1)$ &  &  \\
  &  & $(0,1)$ &  &  &  \\
  &  & $(1,0)$ &  &  &  \\ \hline
\end{tabular}
\end{center}
\caption{Labeling for the 12 small cases which remain after the reductions
 {\bf O}, {\bf T}, {\bf F}, and {\bf R}.\, If $n\equiv 0\pmod 4$, then the center gets
 the unique label occurring fewer in the list than the other elements of $\vi$, and
  if $n\equiv 2\pmod 4$, then it has three options for its label.
 In $(2,2,0,0)$ the center vertex gets the label $(0,0)$; this is an
  exceptional case where only three labels can be used on the non-centers
  and the fourth element of $\vi$ can occur only on the center (cf.\
  Proposition~\ref{p:match}).}\label{tab:5}
}
\end{table}


\subsection{Matchings}
   \label{s:match}

\bigskip


\def \mmm {{\mathcal{M}}}
\def \pmm {{\mathcal{M}_0}}
\def \bbF {{\bf F\boldmath{$'$}}}
\def \bbR {{\bf R\boldmath{$'$}}}
\def \bbT {{\bf T\boldmath{$'$}}}
\def \bno {{\boldmath{$\times$}}}

Recall that a matching (also called packing) in a hypergraph is
 a collection of mutually disjoint edges.
  We now consider hypergraphs whose entire edge set is a matching.
Contrary to the previous parts of the paper, in this particular section
 we allow singleton edges (edges consisting of just one vertex),
  and either exclude or allow isolated vertices.
Let us denote by $\mmm$ the class of hypergraphs with maximum degree 1, i.e.\
 hypergraphs whose edge set is a matching, possibly together with one or more
 vertices of degree 0.
More restrictively let $\pmm\subset \mmm$ denote the subclass
 consisting of the  1-regular hypergraphs, the subscript indicating that
 the number of 0-degree vertices is zero.

Despite that the removal of the center from a star does not change the relative
 value of edge sums --- equal edge sums remain equal, distinct ones remain
 distinct --- this operation is not invariant with respect to $\vi$-cordiality.
This fact, supported by an infinite family of examples,
 is expressed in the following proposition as opposed to Theorem \ref{t:star}.


\begin{proposition}   \label{p:match}
 If $H\in\pmm$ is a hypergraph consisting of mutually disjoint edges, such that
  both $|H|$ and $\|H\|$ are even, moreover
   $|H|\not\equiv \|H\| \pmod 4$, then $H$ is not $\vi$-cordial.
\end{proposition}
\textit{Proof.}\\
Let $E(H)=\{e_1,\dots,e_m\}$ and
 $V(H)=\{v_1,\dots,v_n\}=e_1\cup\cdots\cup e_m$, where the edges
 $e_1,\dots,e_m$ are mutually disjoint.
Consider any vertex labeling $c:V(H)\to\vi$ and its induced
 edge labeling $c^*:E(H)\to\vi$.

Assume that the labeling is $\vi$-cordial, i.e.\ $c$ is
 $\vi$-friendly on $V(H)$ and $c^*$ is $\vi$-friendly on $E(H)$.
Since each vertex belongs to precisely one edge, the sum $S$ of all labels
 satisfies
  $$
    S = \sum_{i=1}^{n} c(v_i) = \sum_{j=1}^{m} c^*(e_j) .
  $$
Now the conditions on $|H|$ and $\|H\|$ imply that precisely one
 of the order and size is a multiple of 4, the other is congruent to $2 \pmod 4$.
For the multiple of 4, every element of $\vi$ occurs the same
 number of times as a vertex label or as an edge label, thus
  $$
    S = (0,0) .
  $$
On the other hand, in the ``\,$2 \pmod 4$\,'' set precisely two
 elemens of $\vi$ occur one fewer times than the other two elements.
Since the overall sum of labels should also be $S=(0,0)$, it follows that
 the sum of two distinct $a,b\in\vi$ should be zero,
 which is impossible.
~\qed

\bigskip

It turns out that this proposition characterizes the exceptions,
 apart from which all matchings are $\vi$-cordial.

\begin{theorem}
 Let $H$ be a matching, where 1-element edges are also allowed.
  \begin{itemize}
    \item[$(i)$] If $H\in\pmm$, then $H$ is $\vi$-cordial if and only if
     $H$ does not satisfy the conditions of Proposition \ref{p:match};
     i.e., if either at least one of $|H|$ and $\|H\|$ is odd,
     or both are even and $|H|\equiv \|H\| \pmod 4$
    \item[$(ii)$] If $H\in\mmm\setminus\pmm$, then $H$ is $\vi$-cordial.
  \end{itemize}
\end{theorem}
\textit{Proof.}\\
 Let $H=(V,E)$, with $n$ vertices and $m$ edges, say
  $E=\{e_1,\dots,e_m\}$.
 The argument mostly applies the ideas of the proof of Theorem \ref{t:star},
  keeping in mind that now $e_m^-=e_m$ holds in any indexing order of the edges.
If $H\in\pmm$, then $H$ can be extended to a star $H^+$ by inserting a center vertex,
 say $x$ ($x\notin V$), and enlarging each edge $e_i$ to $e_i^+:=e_i\cup\{x\}$.
  We already know that $H^+$ has a $\vi$-cordial labeling $c^+$.
If $H$ itself is not $\vi$-cordial, then it must be the case that the label
 of the center occurs one fewer than the most frequent vertex label;
 otherwise we would simply forget about the center and its label.
We are going to prove that this situation can be avoided, unless the
 conditions of Proposition \ref{p:match} hold.

In the same way as in the proof of Theorem \ref{t:star}, one can verify
 that the following reductions are feasible inside the class $\mmm$.
For easier comparison we keep the sequence of properties in same order.
  \begin{enumerate}
    \item If $|e_i|\ge 5$ for some $1\le i\le m$, then we can reduce
     $n$ to $n-4$ by assigning each element of $\vi$ to one vertex of $e_i$,
     while the status of the conditions with respect to $|H|$ and $\|H\|$
     remain unchanged.
      This eliminates all edges larger than 4.
    \item If $|e_1|=|e_2|=|e_3|=|e_4|$, then we can apply the labeling
     scheme given in Table~\ref{tab:2} inside these four edges.
      Then $n$ decreases by a multiple of 4, and $m$ decreases by exactly 4.
      Hence again the conditions with respect to $|H|$ and $\|H\|$
     remain unchanged.
    \item If $|e_1|+|e_2|+|e_3|+|e_4|$ equals 8 or 12, then we can apply the labeling
     scheme given in Table~\ref{tab:3} inside these four edges.
      More explicitly, this step is applicable whenever
       the edges can be indexed in such a way that the sequence
      $(|e_1|,|e_2|,|e_3|,|e_4|)$ is one of
      $(1,1,2,4)$, $(1,1,3,3)$, $(1,2,2,3)$, $(1,3,4,4)$, $(2,2,4,4)$, $(2,3,3,4)$.
      Then again $n$ decreases by a multiple of 4, and $m$ decreases by exactly 4.
      Hence the conditions with respect to $|H|$ and $\|H\|$
     remain unchanged.
    \item If all edges are singletons, or if $H$ has only one edge,
     an obvious labeling verifies that $H$ is $\vi$-cordial.
     Note that in these cases the conditions of Proposition \ref{p:match}
      do not hold because here we have either $|H|=\|H\|$ or $|H|=1$.
    \item If $|e_1|=4$ and $|e_2|>1$, then $(0,0)$ can be assigned to a
     vertex of $e_2$, and the other three elements of $\vi$ to vertices of $e_1$,
     hence inserting partial sums zero in both and reducing$n$ to $n-4$,
     while keeping $m$ unchanged.
     Since $n$ and $m$ do not change modulo 4, the status of the conditions
      on $|H|$ and $\|H\|$ remains the same.
  \end{enumerate}
Steps 1--3 of this list are analogous to (1)--(3) in the proof of
 Theorem \ref{t:star}, while the parts 4 and 5 correspond to the
 reductions {\bf O} and {\bf F}, respectively.

Hence only some of those 49 cases remain to be considered which are marked with
 {\bf T} or {\bf R} or {\bf *} in Table \ref{tab:4}.
For the case of matchings they are summarized in Table \ref{tab:6}.
Among them there are 14 further ones which are reducible by step 5;
 we indicate them with \bbF.
This leaves 35 cases, among which there are 6 satisfying the congruence
 conditions of Proposition \ref{p:match} and hence we know that they are
 not $\vi$-cordial.
These are marked with \bno.

\renewcommand{\arraystretch}{1.5}
\begin{table}[htp]
{
\small
	\begin{center}%
\begin{tabular}
[c]{|c|c|c||c|c|c||c|c|c||c|c|c|
}\hline
 $(0,0,0,2)$ &  $2,8$ & \bbF   & $(0,2,2,0)$ &  $4,10$ & \bno  & $(1,0,3,1)$ &  $5,14$ & \bbF   & $(2,0,0,3)$ &  $5,14$ & \bbF   \\
 $(0,0,0,3)$ &  $3,12$ & \bbF   & $(0,2,3,0)$ &  $5,13$ & \bbT   & $(1,1,0,0)$ &  $2,3$ & \bbR   & $(2,0,1,0)$ &  $3,5$ & \bbT   \\
 $(0,0,2,0)$ &  $2,6$ & \bbR   & $(0,3,0,0)$ &  $3,6$ & {\bf **}   & $(1,1,0,3)$ &  $5,15$ & \bbF   & $(2,2,0,0)$ &  $4,6$ & \bno  \\
 $(0,0,2,3)$ &  $5,18$ & \bbF   & $(0,3,3,0)$ &  $6,15$ & \bbR   & $(1,1,1,0)$ &  $3,6$ & \bbR   & $(2,3,0,0)$ &  $5,8$ & \bbT  \\
 $(0,0,3,0)$ &  $3,9$ & \bbT   & $(0,3,1,0)$ &  $4,9$ & \bbT   & $(1,1,2,0)$ &  $4,9$ & \bbT   & $(3,0,0,1)$ &  $4,7$ & \bbR  \\
 $(0,0,3,2)$ &  $5,17$ & \bbF   & $(0,3,1,1)$ &  $5,13$ & \bbF   & $(1,1,3,0)$ &  $5,12$ & \bbT   & $(3,0,0,2)$ &  $5,11$ & \bbF  \\
 $(0,1,1,0)$ &  $2,5$ & \bbT  & $(0,3,2,0)$ &  $5,12$ & \bbT   & $(1,2,0,0)$ &  $3,5$ & \bbT   & $(3,0,1,0)$ &  $4,6$ & \bno  \\
 $(0,1,1,3)$ &  $5,17$ & \bbF   & $(1,0,0,1)$ &  $2,5$ & \bbT   & $(1,3,0,0)$ &  $4,7$ & {\bf **}   & $(3,0,1,1)$ &  $5,10$ & \bbF  \\
 $(0,1,2,0)$ &  $3,8$ & {\bf **}   & $(1,0,1,0)$ &  $2,4$ & \bno   & $(1,3,0,1)$ &  $5,11$ & \bbF   & $(3,1,0,0)$ &  $4,5$ & \bbT  \\
 $(0,1,3,0)$ &  $4,11$ & \bbR  & $(1,0,1,1)$ &  $3,8$ & \bbF   & $(2,1,0,0)$ &  $3,4$ & {\bf **}  & $(3,1,1,0)$ &  $5,8$ & \bbT  \\
 $(0,2,0,0)$ &  $2,4$ & \bno   & $(1,0,2,0)$ &  $3,7$ & {\bf **}   & $(2,1,1,0)$ &  $4,7$ & \bbR  & $(3,2,0,0)$ &  $5,7$ & \bbT  \\
 $(0,2,0,1)$ &  $3,8$ & \bbF  & $(1,0,3,0)$ &  $4,10$ & \bno   & $(2,0,0,1)$ &  $3,6$ & \bbR   & $(3,3,0,0)$ &  $6,9$ & \bbT  \\
 $(0,2,1,0)$ &  $3,7$ & \bbR & & & & & & & & &  \\
 \hline
\end{tabular}
\end{center}
\caption{The 4-tuples $(m_1,m_2,m_3,m_4)$ not eliminated by steps 1--6,
 the pairs $m,n$, and a way how they can be settled.}\label{tab:6}
}
\end{table}

Note that in the current situation we have $n=f_1+f_2+f_3+f_4$,
 without the $+1$ term; this is the reason why the pairs $m,n$ differ by 1
 when compared in Tables \ref{tab:4} and \ref{tab:6}.
Now a natural analogue of {\bf T} is the following reduction, which
 necessarily is slightly more restrictive.
  \begin{itemize}
    \item[\bbT] --- Trivial reduction applies if
     we have $n\equiv 1\pmod 4$ or $m\equiv 1\pmod 4$
     or both, and $H$ contains an edge whose deletion (also deleting
     its vertices) does not lead to a case marked with \bno.
  \end{itemize}
The reason is that the last vertex can get any label when
 we have a completely balanced labeling on $n-1$ vertices, hence
 the needed label on the last edge can surely be generated; or,
 the last edge can get any label, hence any $\vi$-friendly extension
 of a $\vi$-cordial labeling of the hypergraph with $m-1$ edges
 will do the job.
This operation settles 15 further cases.

\vbox{\noindent As a further simplification, Theorem \ref{l:ext} leads
 to the following reduction.
  \begin{itemize}
    \item[\bbR] --- If there is a non-singleton edge $e_i$ such that
     $H-e_i$ is a matching not marked with \bno, then the following
     conditions are sufficient for reduction:
       $n\equiv 2\pmod 4$ unless $m\equiv 0\pmod 4$, or
       $n\equiv 3\pmod 4$ unless $|e_i|=2$ and $m\equiv 0\pmod 4$.
  \end{itemize}
This eliminates 9 further cases.
}

  \begin{itemize}
    \item[{\bf **}] --- There are 5 cases not covered by the previous
     considerations; Table \ref{tab:7} exhibits a $\vi$-cordial labeling
     for each of them.
  \end{itemize}
This completes the proof of the theorem.
~\qed

\bigskip

\renewcommand{\arraystretch}{1.2}
\begin{table}[htp]
{
  \footnotesize
	\begin{center}%
\begin{tabular}
[c]{|c|c||c|c|c|}\hline
$(m_1,m_2,m_3,m_4)$ & $m,n$ & $f_i = 1$ & $f_i = 2$ & $f_i = 3$ \\ \hline\hline
 $(0,1,2,0)$ & $3,8$ &  & $(0,0)$, $(1,1)$ & $(0,0)$, $(0,1)$, $(1,1)$ \\
  &  &  &  & $(0,1)$, $(1,0)$, $(1,0)$ \\ \hline
 $(0,3,0,0)$ & $3,6$ &  & $(0,0)$, $(0,0)$ &  \\
  &  &  & $(0,1)$, $(1,1)$ &  \\
  &  &  & $(1,0)$, $(1,1)$ &  \\ \hline
 $(1,0,2,0)$ & $3,7$ & $(0,1)$ &  & $(0,0)$, $(0,1)$, $(1,0)$ \\
  &  &  &  & $(1,0)$, $(1,1)$, $(1,1)$ \\ \hline
 $(1,3,0,0)$ & $4,7$ & $(0,0)$ & $(0,1)$, $(1,0)$ & \\
  &  &  & $(0,1)$, $(1,1)$ &  \\
  &  &  & $(1,0)$, $(1,1)$ &  \\ \hline
 $(2,1,0,0)$ & $3,4$ & $(0,1)$ & $(0,0)$, $(1,1)$ &  \\
  &  & $(1,0)$ & &  \\ \hline
\end{tabular}
\end{center}
\caption{Labeling for the 5 final cases of matchings.
 (Edges of size 4 do not occur.)}\label{tab:7}
}
\end{table}

\subsection{Paths}

Inside the class of path hypergraphs we define a \textit{hyperpath}
 as a path in which all edges have size at least~3.
The main result of this section is that every hyperpath is \vcor.
Before proving this, we exhibit an infinite family of paths
 which are not \vcor, hence showing that edges of size 2 create
 more problems than the sporadic examples $P_4$ and $P_5$ themselves.
The complete characterization of $\vi$-cordial paths remains open.


\begin{proposition}   \label{p:path}
 If $H$ is a path with three edges $e_1,e_2,e_3$, such that
  $e_2$ is the middle edge having size $|e_2|=2$, moreover
  $|H|\equiv 0 \pmod 4$, then $H$ is not $\vi$-cordial.
\end{proposition}
\textit{Proof.}\\
Let $V(H)=\{v_1,\dots,v_n\}$, and consider any $\vi$-friendly
 vertex labeling $c:V(H)\to\vi$ with the corresponding induced
 edge labeling $c^*:E(H)\to\vi$.
Since $e_1\cup e_3 = V(H)$ and $|H|$ is a multiple of 4, we now have
  $$
    c^*(e_1) + c^*(e_3) = \sum_{i=1}^{n} c(v_i) = (0,0) .
  $$
This implies $c^*(e_1) = c^*(e_3)$, hence the labeling cannot be
 $\vi$-cordial.
~\qed

\begin{theorem}     \label{t:hyppath}
 Every hyperpath is \vcor.
\end{theorem}
\textit{Proof.}\\
Consider a hyperpath $H=(V,E)$, with $E=\{e_1,\dots,e_m\}$.
 We apply induction on the number $m$ of edges, from $m-4$ to $m$.
The base of induction will be $m=1,2,3$;
 and a special interpretation will be given to the case $m=0$
 to make it possible that the inductive step works for $m=4$,
 hence avoiding the need to verify the assertion separately for the
 many different paths with four edges.

Inside this proof, we simplify the notation to denote the three
 elements of $\vi\setminus \{(0,0)\}$ by $a,b,c$ and write 0 for $(0,0)$.

\medskip

\noindent
\textbf{Case \boldmath$m=1$. }
 Every $\vi$-friendly labeling is \vcor.

\medskip

\noindent
\textbf{Case \boldmath$m=2$. }
 Sequentially creating a $\vi$-friendly labeling,
 for the last vertex we still have at least two choices
 --- which ensure that the sums on $e_1$ and $e_2$ can be made different ---
  unless $n\equiv 0$ (mod~4).
In this exceptional case, however, the sum over the vertex set
 is equal to $0\in\vi$.
Then we assign a nonzero element $b$ to the vertex $e_1\cap e_2$;
 this guarantees that the two sums differ, because the sum over $e_1$
 plus the sum over $e_2$ is equal to $b$.

\medskip

\noindent
\textbf{Case \boldmath$m=3$. }
Let us start with the periodic labeling $0,a,b,c,0,a,b,c,\dots$
 along the vertices of the path, and see whether the sums $s_1,s_2,s_3$ on
  $e_1,e_2,e_3$ are distinct or not.
If some equalities occur, we eliminate them in two steps as follows.

First, to eliminate $s_1=s_3$ if it occurs, we switch the label
 between vertex $e_1\cap e_2$ and its successor (which is only
 in $e_2$, not in $e_1\cup e_3$, because $|e_2|\ge 3$).
This keeps $s_2$ (and also $s_3$) unchanged, but modifies the sum over $e_1$
 to a new updated value of $s_1$, which is then different from $s_3$.

Second, to maintain $s_1\ne s_3$ and eliminate $s_1=s_2$ or $s_2=s_3$
 if it holds after the first step, we switch the label between vertex
 $e_2\cap e_3$ and one of its next two successors.
  (Recall that $|e_3|\ge 3$ holds, hence $|e_3^-|\ge 2$.)
These are two possibilities, each keeping $s_3$ (and also $s_1$)
 unchanged, but offering two new values for an updated $s_2$.
At least one of the two will be different from both $s_1$ and $s_3$,
 hence satisfying the requirement.
(After any of the two switches the original equality $s_1=s_2$ or $s_2=s_3$
automatically disappears, we only have to ensure that a new equality with the
other end will not arise.)

\medskip

\noindent
\textbf{Inductive step from \boldmath$m-4$ to $m$. }
 Instead of dealing with the last four edges, we omit the first two
  and last two edges from the hyperpath $e_1,\dots,e_m$.
Hence let $H'$ be the hyperpath with vertex set $X' = \bigcup_{j=3}^{m-2} e_j$
 and edge set $E' = \{ e_j \mid 3 \le j \le m-2 \}$, with
 $|X'|=n'$ and $|E'|=m'=m-4$.
By the induction hypothesis there exists a \vcor\ labeling $(c',c'^*)$ on
 $(X',E')$.
Our goal is to assign $n-n'$ labels to the vertices of $V\setminus X'$
 and generate four distinct sums on $e_1,e_2,e_{m-1},e_m$.
The $n-n'$ labels have to be selected from a balanced multiset $S'$ of
 $4\cdot\lceil n/4 \rceil - n'$ elements over $\vi$; namely, starting with
 $\lceil n/4 \rceil$ copies of $\vi$ we delete the elements which have been
 assigned to $X'$, and from the remaining multiset we need to select
 $n-n'$ labels properly.
(Note that the multiplicities of  any two elements in $S'$ differ by at most 1,
 because $c'$ is $\vi$-friendly by assumption, hence what remains after
  omitting them from $\lceil n/4 \rceil$ times $\vi$ is also balanced.)
We can assume without loss of generality that $n\equiv 0$ (mod~4), because
 any other case would give us some flexibility in selecting the set of labels,
 whereas in this case the multiset of labels to be used is determined.

Assume that the vertices $e_2\cap X'$ and $e_{m-1}\cap X'$ are labeled
 with $x$ and $y$, respectively, and that the sum of all labels over $X'$
 is $z$.
  (Some or all of $x,y,z$ may coincide.)
Then the label $x'$ of $e_1\cap e_2$ and $y'$ of $e_{m-1}\cap e_m$
 must satisfy
  \begin{equation}   \label{e:xy}
    x' + y' = x + y + z .
  \end{equation}
Indeed, since $n\equiv 0$ (mod~4), the sum $z$ of labels inside $X'$ is equal to
 the sum outside $X'$, moreover the sums over the edges $e_1,e_2,e_{m-1},e_m$
 is equal to $0\in\vi$, in which $x+y+x'+y'$ is counted additionally to $z$.

We proceed in three steps, after which a \vcor\ labeling will be obtained.

\smallskip

\textit{Step 1.} Determine $x',y'$.

We choose $x'$ and $y'$ in such a way that one of them is an element which is
 most frequent in $S'$, moreover the remaining multiset $S'\setminus\{x',y'\}$
 still contains at least one occurrence of $0$.
  We argue that this can always be done.
Indeed, the condition on edge sizes implies $|S'|\ge 8$.
Assume first that equality $|S'|=8$ holds; then
 each element occurs precisely twice in $S'$.
If equation~(\ref{e:xy}) requires $x'=y'$ (that is, if $x+y+z=0$), then we can use
 any of the three labels different from $0$ for $x'$.
On the other hand if $x'\ne y'$, the required sum $x'+y'$ can be formed in two ways,
 each of them leaving two elements of $\vi$ with multiplicity 2 and two with 1
 in $S'\setminus\{x',y'\}$, hence either choice is feasible.
Finally if $|S'|>8$, the most frequent element occurs at least three times.
We choose it for $x'$, and assign $x+y+z-x'$ for $y'$.
This is feasible because all elements have multiplicity
 at least 2 in $S'\setminus\{x'\}$.

\smallskip

\textit{Step 2.} Distribute all but 6 labels from $S'\setminus\{x',y'\}$.

If $|S'|=8$, there is nothing to do in this step, the remaining multiset is
 $$
   0,0,a,a,b,b \quad \mbox{\rm or} \quad 0,a,a,b,b,c
    \quad \mbox{\rm or} \quad 0,0,a,a,b,c
 $$
  whose sum is
 $$
   0 \quad \mbox{\rm or} \quad c
    \quad \mbox{\rm or} \quad a
 $$
  respectively.
If $|S'|>8$, we distribute $|S'|-8$ elements from $S'\setminus\{x',y'\}$
 almost arbitrarily, but
 in such a way that the following conditions are met:
 \begin{itemize}
   \item either $0,a,a,b,b,c$ or $0,0,a,a,b,c$ remains;
   \item  $e_2$ and $e_{m-1}$ have just one unlabeled vertex;
   \item each of $e_1$ and $e_m$ has two unlabeled vertices.
 \end{itemize}
After this, let us denote the current sums of labels in
 $e_1,e_2,e_{m-1},e_m$ by $s_1,s_2,s_{m-1},s_m$, respectively.
From these four partial sums we shall have to create four distinct
 final sums by properly distributing the remaining six labels.
From this point of view $(s_1,s_2,s_{m-1},s_m)$ and
 $(a+s_1,a+s_2,a+s_{m-1},a+s_m)$ are equivalent.
For this reason we may assume without loss of generality that
 0 is most frequent among $s_1,s_2,s_{m-1},s_m$.
Hence, apart from the permutation of subscripts, only the
 following five types of 4-tuples are relevant for $(s_1,s_2,s_{m-1},s_m)$.
 \begin{enumerate}
  \item --- $(0,0,0,0)$, sum = 0
  \item --- $(0,0,0,a')$, sum = $a'\neq 0$
  \item --- $(0,0,a',a')$, sum = 0
  \item --- $(0,0,a',b')$, sum = $a'+b'\neq 0$
  \item --- $(0,a',b',c')$, sum = $a'+b'+c'$ = 0
 \end{enumerate}
Here we use prime notation to mean that $a',b'$ may be other than
 $a,b$ in the remaining 6-element multiset; but different
 primed letters mean different elements.
Observe, however, that $s_1+s_2+s_{m-1}+s_m$ must be equal to the sum
 of the six elements in the multiset, because the total sum over the four
 edges will eventually be zero; this is implied by the choice
 of $x'$ and $y'$.
This fact yields, in particular, that not all 4-tuples fit together with
 all 6-tuples.
Namely, $0,0,a,a,b,b$ is compatible with the cases $1,3,5$ while
 $0,a,a,b,b,c$ and $0,0,a,a,b,c$ admit the cases $2,4$.

\smallskip

\textit{Step 3.} Complete the labeling on $e_1,e_2,e_{m-1},e_m$.

 This step is a little time consuming, but easy.
The selection rules described above already imply that if three edges
 have mutually distinct final sums, then the fourth edge has the
 missing value for its sum.
To achieve this, we systematically enumerate the 4-tuples listed in
 1--5 above with their compatible 6-tuples of labels, and --- up to
  symmetry --- the possible positions of $0,a',b'$ and the
  elements that can play the role of $0$, $a'$, and $b'$.
Tables \ref{tab:p-1} and \ref{tab:p-2} exhibit a suitable way of
 extending $c'$ to a \vcor\ labeling of the entire path.

\medskip

\noindent
\textbf{Case \boldmath$m=4$. }
  Let us
 artificially introduce the 0-path as a single vertex with no edges.
  It is of course \vcor, any label $x$ can be assigned to the vertex.
Now, for $m=4$ we identify the vertex with $e_2\cap e_3$, and
 apply the inductive step above as described for the case $x=y$.
This completes the proof of the theorem.~\qed

\bigskip

\renewcommand{\arraystretch}{1.2}
\begin{table}[htp]
{\footnotesize
	\begin{center}%
\begin{tabular}
[c]{|c|c|c||c|c|c|c|}\hline
$(s_1,s_2,s_{m-1},s_m)$ & 6-tuple & $a'=$ & $e_1$ & $e_2$ & $e_{m-1}$ & $e_m$ \\ \hline\hline
 $(0,0,0,0)$ & $0,0,a,a,b,b$ & $-$ & $0,a \to a$ & $0 \to 0$ & $b \to b$ & $a,b \to c$ \\
 \hline
 $(0,0,0,a')$ & $0,a,a,b,b,c$ & $c$ & $a,b \to c$ & $0 \to 0$ & $b \to b$ & $a,c \to a$ \\
  & $0,0,a,a,b,c$ & $a$ & $0,a \to a$ & $0 \to 0$ & $b \to b$ & $a,c \to c$ \\
 $(0,0,a',0)$ & $0,a,a,b,b,c$ & $c$ & $a,b \to c$ & $0 \to 0$ & $b \to a$ & $a,c \to b$ \\
  & $0,0,a,a,b,c$ & $a$ & $0,a \to a$ & $0 \to 0$ & $b \to c$ & $a,c \to b$ \\
 \hline
 $(0,0,a',a')$ & $0,0,a,a,b,b$ & $a$ & $0,a \to a$ & $0 \to 0$ & $b \to c$ & $a,b \to b$ \\
  & & $c$ & $0,b \to b$ & $0 \to 0$ & $b \to a$ & $a,a \to c$ \\
 $(0,a',0,a')$ & & $a$ & $a,b \to c$ & $0 \to a$ & $b \to b$ & $0,a \to 0$ \\
  & & $c$ & $0,a \to a$ & $0 \to c$ & $b \to b$ & $a,b \to 0$ \\
 $(0,a',a',0)$ & & $a$ & $0,b \to b$ & $0 \to a$ & $b \to c$ & $a,a \to 0$ \\
  & & $c$ & $0,b \to b$ & $0 \to c$ & $b \to a$ & $a,a \to 0$ \\
 $(a',0,0,a')$ & & $a$ & $0,b \to c$ & $0 \to 0$ & $b \to b$ & $a,a \to a$ \\
  & & $c$ & $0,b \to a$ & $0 \to 0$ & $b \to b$ & $a,a \to c$ \\
 \hline
\end{tabular}
\end{center}
\caption{Labels inserted into $e_1,e_2,e_{m-1},e_m$
 starting from at most two distinct sums, and the final sum of labels
  inside $e_i$.}\label{tab:p-1}
}
\end{table}

\renewcommand{\arraystretch}{1.2}
\begin{table}[htp]
{\footnotesize
	\begin{center}%
\begin{tabular}
[c]{|c|c|c||c|c|c|c|}\hline
$(s_1,s_2,s_{m-1},s_m)$ & 6-tuple & $a',b'=$ & $e_1$ & $e_2$ & $e_{m-1}$ & $e_m$ \\ \hline\hline
 $(0,0,a',b')$ & $0,a,a,b,b,c$ & $a,b$ & $a,a \to 0$ & $b \to b$ & $b \to c$ & $0,c \to a$ \\
  & $0,0,a,a,b,c$ & $b,c$ & $a,a \to 0$ & $b \to b$ & $c \to a$ & $0,0 \to c$ \\
 $(0,a',0,b')$ & $0,a,a,b,b,c$ & $a,b$ & $a,a \to 0$ & $b \to c$ & $b \to b$ & $0,c \to a$ \\
  & $0,0,a,a,b,c$ & $b,c$ & $a,a \to 0$ & $0 \to b$ & $c \to c$ & $0,b \to a$ \\
 $(0,a',b',0)$ & $0,a,a,b,b,c$ & $a,b$ & $a,a \to 0$ & $b \to c$ & $c \to a$ & $0,b \to b$ \\
  & $0,0,a,a,b,c$ & $b,c$ & $a,a \to 0$ & $0 \to b$ & $0 \to c$ & $b,c \to a$ \\
 $(a',0,0,b')$ & $0,a,a,b,b,c$ & $a,b$ & $b,c \to 0$ & $a \to a$ & $b \to b$ & $0,a \to c$ \\
  & $0,0,a,a,b,c$ & $b,c$ & $0,b \to 0$ & $a \to a$ & $c \to c$ & $0,a \to b$ \\
 \hline
 $(0,a',b',c')$ & $0,0,a,a,b,b$ & $a,b$ & $a,a \to 0$ & $0 \to a$ & $0 \to b$ & $b,b \to c$ \\
  & & $b,c$ & $a,a \to 0$ & $0 \to b$ & $0 \to c$ & $b,b \to a$ \\
 $(a',0,b',c')$ & & $a,b$ & $0,a \to 0$ & $a \to a$ & $0 \to b$ & $b,b \to c$ \\
  & & $b,c$ & $0,b \to 0$ & $b \to b$ & $0 \to c$ & $a,a \to a$ \\
 \hline
\end{tabular}
\end{center}
\caption{Labels inserted into $e_1,e_2,e_{m-1},e_m$
 starting from 3 or 4 distinct sums, and the final sum of labels
  inside $e_i$.}\label{tab:p-2}
}
\end{table}

\section{Conclusions}

We finish the paper with some simple open problems.

\begin{conjecture}
Let $T=(V,E)$ be a hypertree.
 If $|e|\geq3$ for every $e\in E(T)$, then $T$ is $\vi$-cordial.
\end{conjecture}


\begin{problem}
 Characterize the class of hypergraphs which are cycle-free and $\vi$-cordial.
\end{problem}

\begin{problem}
 Give necessary and sufficient conditions for $\vi$-cordial path hypergraphs.
\end{problem}

\section{Acknowledgments}
Research of Sylwia Cichacz and Agnieszka G\"orlich was partially supported by the Faculty of Applied Mathematics AGH UST statutory tasks within subsidy of Ministry of Science and Higher Education. Research of Zsolt Tuza was supported by the National Research, Development and
 Innovation Office -- NKFIH under the grant SNN 116095.



\end{document}